\documentclass[11pt]{article}
\usepackage[T1]{fontenc}
\usepackage{lmodern}
\usepackage{amsmath,amssymb,amsthm}
\usepackage[a4paper,margin=28mm]{geometry}
\usepackage[shortlabels]{enumitem}
\usepackage{array,graphicx,needspace,longtable}
\usepackage{xcolor} 
\usepackage[colorlinks=true,citecolor=blue,linkcolor=blue,urlcolor=blue]{hyperref}

\makeatother
\newtheorem{theorem}{Theorem}[section]
\newtheorem{proposition}[theorem]{Proposition}
\newtheorem{lemma}[theorem]{Lemma}
\newtheorem{corollary}[theorem]{Corollary}
\theoremstyle{definition}
\newtheorem{definition}[theorem]{Definition}
\newtheorem{example}[theorem]{Example}
\newtheorem{remark}[theorem]{Remark}
\DeclareMathOperator{\spec}{spec}

\DeclareMathOperator{\diag}{diag}

\newcommand{\BH}{\mathrm{BH}}
\newcommand{\one}{\mathbf{1}}
\newcommand{\ii}{\mathrm{i}}
\newcommand{\R}{\mathcal{R}}
\newcommand{\C}{\mathcal{C}}

\title{From Butson Hadamard switching to graph switching and back}
\author{Aida Abiad\thanks{\texttt{a.abiad.monge@tue.nl}, Department of Mathematics and Computer Science, Eindhoven University of Technology, The Netherlands; Department of Mathematics and Data Science, Vrije Universiteit Brussel, Belgium.}
\and Ronan Egan\thanks{\texttt{ronan.egan@dcu.ie}, School of Mathematical Sciences, Dublin City University, Ireland.}
\and Stef Mohnen\thanks{\texttt{s.p.j.mohnen@student.tue.nl}, Department of Mathematics and Computer Science, Eindhoven University of Technology, The Netherlands.}}
\date{}

\begin{document}
\maketitle

\begin{abstract}
Switching is a local operation that preserves specified properties of a combinatorial object. It has been used to derive constructions of inequivalent Hadamard matrices and nonisomorphic cospectral graphs. In this paper we investigate the relation between these constructions for real and Butson Hadamard matrices. We show that, while equivalent matrices have isomorphic associated graphs, the converse can fail for Butson matrices. We prove that a character decomposition determines the graph spectrum from the singular values of entrywise powers. For prime phase, we also observe that the spectrum depends only on the order and phase; for composite phase, a switching preserving the Butson property can change it. We show that Orrick's closed-quadruple and Hall-set switchings correspond to Wang–Qiu–Hu switching, up to graph isomorphism. Going in the other direction, we characterise when the Godsil–McKay and Wang–Qiu–Hu reflections preserve the Butson property and give phase-multiplicity conditions for the corresponding graph switchings. Finally we also relate block switching to cyclic graph switching under a phase-balance condition, automatic at prime phase. These results provide a basis for transferring tools and techniques between matrix switching and graph switching.
\end{abstract}

\noindent\textbf{2020 Mathematics Subject Classification.} Primary 05B20; Secondary 05C50, 15B34.

\smallskip
\noindent\textbf{Keywords.} Butson Hadamard matrix; cospectral graphs; Godsil--McKay switching; Wang--Qiu--Hu switching; closed quadruple; Hall set.

\section{Introduction}\label{section:introduction}

Switching is a local operation that changes a combinatorial object while preserving specified properties. Switching methods have been used to construct inequivalent Hadamard matrices and nonisomorphic cospectral graphs, that is, graphs with the same adjacency spectrum. We study the relation between these constructions. Our aim is to determine when matrix switching induces switching of an associated graph and how graph-switching methods yield operations on matrices.

Butson Hadamard matrices (sometimes called generalised Hadamard matrices) generalise real Hadamard matrices by allowing entries in a fixed finite group of roots of unity \cite{butson1962}; see \cite{craigen2006hadamard,deLF} for background. Orrick \cite{orrick2008switching} introduced two switching operations for real Hadamard matrices, based on closed quadruples and Hall sets. Crnkovi\'{c}, Egan and \v{S}vob \cite{GenSwitching} developed a common framework for switching combinatorial objects and extended these operations to Butson and complex Hadamard matrices. Related constructions include Di\c{t}\u{a}'s block method \cite{Dita2004}, which gives Butson matrices when the input entries and parameters lie in a common finite group of roots of unity. Local deformations also arise from Sz\"oll\H{o}si's paired-column construction \cite{Szollosi2008} and Goyeneche's extension using pairs of rows or columns whose coordinatewise quotients are $\pm1$ \cite{Goyeneche2013}.

In spectral graph theory, Godsil--McKay switching (GM-switching) \cite{GMswitching} and Wang--Qiu--Hu switching (WQH-switching) \cite{WQHswitching} provide methods for constructing cospectral graphs. Abiad and Peters \cite{abiad2024switching} related Orrick's operations to GM-switching on the Hadamard graphs introduced by McKay \cite{McKay1979}. In McKay's construction, two real Hadamard matrices are equivalent if and only if their associated graphs are isomorphic. This connection motivates extending the graph-switching approach to Butson matrices.

In this paper, we use a graph construction that we call Butson Hadamard graph (see Definition~\ref{def:BHGraph}) to study equivalence, spectra and switching for Butson Hadamard matrices. Each matrix entry is replaced by a cyclic matching between a row fibre and a column fibre, with loops distinguishing the two parts. The construction includes McKay's real case. We show that equivalent matrices have isomorphic graphs (Proposition~\ref{prop:EquivIso}), but the converse need not hold for Butson matrices. Field conjugates have isomorphic graphs (Proposition~\ref{prop:FAut}), and Example~\ref{ex:isomorphic_inequivalent} gives matrices whose graphs are isomorphic even though the matrices remain inequivalent when field automorphisms and transposition are allowed.

Next we determine the graph spectra by a character decomposition in Section~\ref{section:spectra}. The spectrum is expressed in terms of the singular values of the entrywise powers of the matrix. For a fixed matrix order and a fixed prime phase, we prove that all associated graphs are cospectral (Theorem~\ref{thm:prime-cospectral}). A switching in a Butson Hadamard matrix of prime phase therefore preserves the graph spectrum when $k$ is prime. For composite phase, this conclusion need not hold; indeed, Example~\ref{ex:k4} gives a switching that preserves the Butson property but changes the graph spectrum. Thus the matrix conditions alone do not always suffice to obtain a corresponding graph switching.

We then study when particular switching methods can be transferred between the two settings. For real Hadamard matrices, Orrick's closed-quadruple and Hall-set switchings correspond, up to graph isomorphism, to WQH-switching (Theorems~\ref{thm:ClosedQuadrupleWQH} and~\ref{thm:HallSetWQH}). Going the other way around, and with the goal to obtain operations on Butson matrices from graph switching, we apply the reflections underlying GM- and WQH-switching. We characterise when these reflections preserve the Butson property and determine the phase-multiplicity conditions under which they induce the corresponding graph switchings on the specified partitions. For the paired WQH reflection, these graph conditions are automatic for every admissible matrix switching at prime phase (Corollary~\ref{cor:butson_wqh_prime}). At even phase, the WQH reflection is related to a phased GM reflection by allowed row scalings. The role of phase multiplicities also appears in block switching. Corollary~\ref{cor:cyclic_BH} shows that the block switching of Proposition~\ref{prop:FourierSwitch} induces the cyclic graph switching of Theorem~\ref{thm:cyclic_switching} when each complementary row has uniform phase multiplicities. This condition is automatic at prime phase. These results provide a basis for transferring tools and techniques between matrix switching and graph switching.

This paper is organised as follows. Section~\ref{section:preliminaries} introduces the notation, graph construction and switching methods. Section~\ref{section:equivalence} studies matrix equivalence, graph isomorphism and reconstruction. Section~\ref{section:spectra} determines the spectra of the associated graphs. Section~\ref{section:new_results} establishes the correspondence between Orrick's switchings and suitable WQH-switchings. Section~\ref{section:butson_switching} studies the GM and WQH reflections for Butson matrices, the corresponding graph switchings, and cyclic switching for balanced blocks. Section~\ref{section:concluding_remarks} discusses further questions.

\section{Preliminaries}\label{section:preliminaries}

\subsection{Notation and Butson Hadamard matrices}\label{subsection:notation}

Throughout we write $\mathbb Z_k=\mathbb Z/k\mathbb Z$, $\zeta_k=e^{2\pi\ii/k}$, and $\mu_k=\{1,\zeta_k,\ldots,\zeta_k^{k-1}\}$. An $n \times n$ matrix $H$ with entries in $\mu_{k}$ such that $HH^{\ast} = nI_{n}$ is \emph{Butson Hadamard} of order $n$ and phase $k$, where $H^{\ast}$ denotes the complex conjugate transpose of $H$. The set of all such matrices is denoted by $\mathrm{BH}(n,k)$. An element of $\BH(n,2)$ is a \emph{Hadamard matrix}, and will sometimes be referred to as a real Hadamard matrix to emphasise this case. A complex Hadamard matrix satisfies the same orthogonality condition and has entries of absolute value one, without the requirement that they are roots of unity. Real Hadamard matrices have order $1$, $2$, or a multiple of $4$. We use $k$ as a specified phase; it need not always be the least phase containing all entries of $H$.

Matrix rows and columns are indexed from $0$ to $n-1$, unless a local block indexing is specified. Indices in $\mathbb Z_k$ are always read modulo $k$. We write $I_s$, $J_{s,t}$ and $O_{s,t}$ for the identity, all-one and zero matrices of the indicated sizes, and abbreviate $J_{s,s}$ and $O_{s,s}$ to $J_s$ and $O_s$. The vector $\one_s$ is the column vector of $s$ ones. In a small number of instances, by abuse of notation $\one_X$ will also denote a binary characteristic vector in $\mathbb{R}^{|S|}$ for some subset $X \subseteq S$, and this will be stated or clear from context in those instances. The symbols $\otimes$ and $\oplus$ denote the Kronecker product and the direct sum, respectively. The Hadamard product $A\circ B$ is the entrywise product of matrices of the same size. For a matrix with entries in $\mu_k$, we write
\[
 A^{(\ell)}=(a_{ij}^{\ell})_{i,j}\qquad(\ell\in\mathbb Z)
\]
for its entrywise $\ell^{\rm th}$ power. In particular, $A^{(0)}$ is an all-ones matrix. The parentheses distinguish entrywise powers from ordinary matrix powers. The spectrum $\spec(A)$ is the multiset of eigenvalues of $A$, counted with algebraic multiplicity. Unions of spectra are understood as multiset unions.

Two matrices $H,K\in\BH(n,k)$ are \emph{equivalent} if
\[
 K=M H N^*
\]
for monomial matrices $M,N$ whose nonzero entries belong to $\mu_k$. A monomial matrix has exactly one nonzero entry in each row and each column. Thus equivalence allows row and column permutations and multiplication of individual rows and columns by elements of $\mu_k$. We call $H$ \emph{normalised} if its first row and first column consist entirely of ones. Every Butson Hadamard matrix is equivalent to a normalised one.

For $\ell$ coprime to $k$, the map $\zeta_k\mapsto\zeta_k^\ell$ is a field automorphism of $\mathbb Q(\zeta_k)$ and sends $H$ to another matrix $H^{(\ell)}\in\BH(n,k)$. We distinguish this operation, and transposition, from the equivalence relation just defined. When we allow both in addition to monomial equivalence, we refer explicitly to \emph{extended equivalence}.

\subsection{Graphs associated with phase matrices}\label{subsection:graph_construction}

All graphs in this paper are finite and undirected, with no multiple edges. We allow loops when explicitly stated. A loop contributes a diagonal entry $1$ to the adjacency matrix and is counted once in the neighbourhood of its vertex. For a graph $G$, a vertex $v$ and a subset $S\subseteq V(G)$, we write
\[
 N_G(v)=\{u\in V(G):u\text{ is adjacent to }v\},\qquad
 d_S(v)=|N_G(v)\cap S|.
\]
We omit the subscript $G$ when the graph is clear. The \emph{spectrum} of a graph is the spectrum of its adjacency matrix. Two graphs are \emph{cospectral} if their spectra agree, and are \emph{cospectral mates} if they are also nonisomorphic.

Let $P_k$ be the permutation matrix indexed by $\mathbb Z_k$ with
\[
 (P_k)_{a,b}=1\quad\Longleftrightarrow\quad b=a+1.
\]
Thus $P_k$ is the circulant permutation matrix with first row $(0,1,0,\ldots,0)$. The map
\[
 \phi_k:\mu_k\longrightarrow\langle P_k\rangle,
 \qquad \phi_k(\zeta_k^a)=P_k^a,
\]
is the regular permutation representation of $\mu_k$. For a matrix $H$ with entries in $\mu_k$, we write
\[
 \phi_k(H)=(\phi_k(h_{ij}))_{i,j}.
\]
This is a block substitution. We do not regard $\phi_k$ as an additive map on complex numbers. Where a monomial matrix is involved, we set $\phi_k(0)=O_k$; products with monomial matrices then respect this substitution because each contributing entry is a single product. We often omit the subscript $k$ when the phase is fixed.

\begin{definition}\label{def:BHGraph}\label{HadamardGraphdef}
Let $H=(\zeta_k^{e_{ij}})\in\mu_k^{n\times n}$. The graph $G_H$ has vertex set $\R\mathbin{\dot\cup}\C$, where
\[
 \R=\{R_{i,a}:0\leq i<n,\ a\in\mathbb Z_k\},\qquad
 \C=\{C_{j,b}:0\leq j<n,\ b\in\mathbb Z_k\}.
\]
Its edges between the two parts are given by
\begin{equation}\label{eq:adjacency_rule}
 R_{i,a}\sim C_{j,b}\quad\Longleftrightarrow\quad b-a=e_{ij}\text{ in }\mathbb Z_k.
\end{equation}
Every vertex in $\R$ has a loop, and there are no other edges. If $H\in\BH(n,k)$, we call $G_H$ the \emph{Butson Hadamard graph} associated with $H$.
\end{definition}

We always use the vertex order
\begin{equation}\label{eq:partition}
 R_{0,0},\ldots,R_{0,k-1},\ldots,R_{n-1,k-1},
 C_{0,0},\ldots,C_{0,k-1},\ldots,C_{n-1,k-1}.
\end{equation}
With this convention, the adjacency matrix is
\begin{equation}\label{eq:AltAD}
 \mathcal M_H=
 \begin{pmatrix}
 I_{nk}&\phi_k(H)\\
 \phi_k(H)^{\top}&O_{nk}
 \end{pmatrix}.
\end{equation}
The graph obtained from $G_{H}$ by deleting the loops is bipartite, with parts $\R$ and $\C$. The loops in $G_H$ distinguish the row part from the column part under graph isomorphism.

For each row $i$ and column $j$, the sets
\[
 R_i=\{R_{i,a}:a\in\mathbb Z_k\},\qquad
 C_j=\{C_{j,b}:b\in\mathbb Z_k\}
\]
are called the \emph{row fibres} and \emph{column fibres}. Each pair $R_i,C_j$ induces a perfect matching. The simultaneous translation $R_{i,a}\mapsto R_{i,a+c}$, $C_{j,b}\mapsto C_{j,b+c}$ preserves adjacency. If $H$ is normalised, then $R_{0,a}\sim C_{j,a}$ and $R_{i,a}\sim C_{0,a}$ for every $i$, $j$ and $a$.

When $k=2$, we recover McKay's graph \cite{McKay1979}: the two copies of row $i$ are $R_{i,0}$ and $R_{i,1}$, the two copies of column $j$ are $C_{j,0}$ and $C_{j,1}$, and a negative entry interchanges the corresponding matching. We use the same vertex notation in the real and complex cases.

\subsection{Construction of cospectral graphs via GM-switching and WQH-switching}\label{subsection:graph_switching}

We recall the two graph-switching methods used below. Their usual matrix proofs also apply under our loop convention.

\begin{theorem}[GM-switching \cite{GMswitching}]\label{thm:butson_gm_classical}
Let $\{X_1,\ldots,X_t,D\}$ be a partition of $V(G)$ such that:
\begin{enumerate}[(i)]
\item for every $i,j$, the number $d_{X_j}(v)$ is independent of the choice of $v\in X_i$;
\item for every $v\in D$ and every $i$, the number $d_{X_i}(v)$ belongs to $\{0,|X_i|/2,|X_i|\}$.
\end{enumerate}
For every $v\in D$ with $d_{X_i}(v)=|X_i|/2$, reverse all adjacencies between $v$ and $X_i$. The resulting graph is cospectral with $G$.
\end{theorem}

The orthogonal matrix underlying this operation acts on each $X_i$ as the reflection
\begin{equation}\label{eq:butson_gm_Q}
 Q_s=\frac{2}{s}J_s-I_s\qquad(s=|X_i|)
\end{equation}
and as the identity on $D$. We have $Q_s^{\top}=Q_s$, $Q_s^2=I_s$ and $Q_s\one_s=\one_s$.

\begin{theorem}[WQH-switching \cite{WQHswitching}]\label{thm:WQH}
Let
\[
 \{D,X_1^+,X_1^-,\ldots,X_t^+,X_t^-\}
\]
be a partition of $V(G)$ such that, for every $i,j$:
\begin{enumerate}[(i)]
\item $|X_i^+|=|X_i^-|$;
\item there is a constant $c_{ij}$ satisfying
\[
 \begin{cases}
 d_{X_j^+}(v)-d_{X_j^-}(v)=c_{ij},&v\in X_i^+,\\
 d_{X_j^-}(v)-d_{X_j^+}(v)=c_{ij},&v\in X_i^-;
 \end{cases}
\]
\item every vertex $v\in D$ is adjacent to all vertices of one of $X_i^+,X_i^-$ and none of the other, or satisfies $d_{X_i^+}(v)=d_{X_i^-}(v)$.
\end{enumerate}
Whenever a vertex of $D$ is adjacent to all of one of $X_i^+,X_i^-$ and none of the other, interchange these two neighbourhoods. The resulting graph is cospectral with $G$.
\end{theorem}

\begin{remark}\label{def:WQH_simplified}
For $t=1$, the conditions concern a pair $(X^+,X^-)$ of equally sized sets. The difference between the number of neighbours in a vertex's own set and in the other set must be constant on $X^+\cup X^-$. Every remaining vertex must have either equal numbers of neighbours in the two sets, or all of one set and none of the other. We use this form for closed-quadruple switching. Hall-set switching requires two pairs of sets.
\end{remark}

For completeness, the loop convention can be checked directly in the matrix formulation. Put $s_i=|X_i^+|$ and $z_i=\one_{X_i^+}-\one_{X_i^-}$, where the latter are characteristic vectors in $\mathbb R^{|V(G)|}$. The WQH conjugating matrix is
\[
 W=I-\sum_{i=1}^t\frac{1}{s_i}z_i z_i^{\top}.
\]
The vectors $z_i$ have disjoint supports and squared norm $2s_i$, so $W$ is orthogonal. Condition~(ii) ensures that the adjacency block on the switching sets is fixed by conjugation. Condition~(iii) makes the remaining blocks either fixed or switched as prescribed. Thus loops are included by counting their diagonal entries once, just as in the stated neighbourhood conditions.

\subsection{Construction of unequivalent Hadamard matrices via Closed quadruples and Hall set switching methods}\label{subsection:orricks_switching}

We recall Orrick's two switching operations \cite{orrick2008switching}. Let $h_i$ denote row $i$ of a real Hadamard matrix $H$.

\begin{definition}\label{3normalized}
A Hadamard matrix is \emph{3-normalised} on three distinct rows $i,j,\ell$ if $h_i\circ h_j\circ h_\ell=\one_n^{\top}$. A set of four distinct rows is a \emph{closed quadruple} if their Hadamard product is $\pm\one_n^{\top}$.
\end{definition}

After equivalence operations, a Hadamard matrix with a closed quadruple can be written as
\begin{equation}\label{closedquad}
 H=\begin{pmatrix}
 \one_m^{\top}&-\one_m^{\top}&-\one_m^{\top}&\one_m^{\top}\\
 \one_m^{\top}&-\one_m^{\top}&\one_m^{\top}&-\one_m^{\top}\\
 \one_m^{\top}&\one_m^{\top}&-\one_m^{\top}&-\one_m^{\top}\\
 \one_m^{\top}&\one_m^{\top}&\one_m^{\top}&\one_m^{\top}\\
 A_1&A_2&A_3&A_4
 \end{pmatrix},\qquad m=n/4.
\end{equation}
The four sets of $m$ columns are called \emph{fields}. Orthogonality implies that every row of each $A_a$ has sum zero. Consequently, if $n>4$, then $m$ is even and $n$ is divisible by $8$.

\begin{definition}
\emph{Switching a closed quadruple} means negating the entries in the four selected rows within one field.
\end{definition}

\begin{theorem}[Orrick \cite{orrick2008switching}]\label{thm:orrick_quad}
Switching a closed quadruple produces a Hadamard matrix. The four choices of field produce equivalent matrices.
\end{theorem}

For $n>4$, a \emph{Hall set of rows} consists of four rows that can be 3-normalised on three of them so that the fourth row has one exceptional sign in each field. Equivalently, their Hadamard product has exactly four entries of one sign and all remaining entries of the other sign. The columns containing these four exceptional entries are the \emph{Hall columns}. For orders divisible by $8$, the Hall columns form a closed quadruple; for orders congruent to $4$ modulo $8$, they form a Hall set \cite{orrick2008switching}. We use \emph{Hall set} for the four rows together with their Hall columns in the latter case.

For $n=4m+4\geq12$, with $m$ even, a Hadamard matrix containing a Hall set is equivalent to one of the form
\begin{equation}\label{Hallsetmatrixblocks}
 H=\begin{pmatrix}
 S_4&F_1&F_2&F_3&F_4\\
 G_1&A_{11}&A_{12}&A_{13}&A_{14}\\
 G_2&A_{21}&A_{22}&A_{23}&A_{24}\\
 G_3&A_{31}&A_{32}&A_{33}&A_{34}\\
 G_4&A_{41}&A_{42}&A_{43}&A_{44}
 \end{pmatrix}.
\end{equation}
Here $S_4=2I_4-J_4$, $F_a=f_a\one_m^{\top}$, where
\[
 f_1=\begin{pmatrix}1\\1\\1\\1\end{pmatrix},\quad
 f_2=\begin{pmatrix}1\\-1\\-1\\1\end{pmatrix},\quad
 f_3=\begin{pmatrix}1\\-1\\1\\-1\end{pmatrix},\quad
 f_4=\begin{pmatrix}-1\\-1\\1\\1\end{pmatrix},
\]
and $G_1=F_1^{\top}$, $G_a=-F_a^{\top}$ for $a=2,3,4$. Each $A_{ab}$ is an $m\times m$ matrix whose row and column sums are $2$ if $a=b$ and $0$ otherwise. 

\begin{definition}\label{switchingHallsetblocks}
\emph{Switching a Hall set} is negating $F_a$ and $G_a$ for one $a\in\{1,2,3,4\}$ in \eqref{Hallsetmatrixblocks}.
\end{definition}

\begin{theorem}[Orrick \cite{orrick2008switching}]\label{thm:orrick_hall}
Switching a Hall set produces a Hadamard matrix. The four choices of the pair $(F_a,G_a)$ produce equivalent matrices.
\end{theorem}

\section{Equivalence and reconstruction}\label{section:equivalence}

We first describe the relation between equivalence of phase matrices and isomorphism of their associated graphs. The graph construction retains monomial equivalence, but it also identifies matrices related by field automorphisms.

\begin{proposition}\label{prop:EquivIso}
If $H,K\in\BH(n,k)$ are equivalent, then $G_H$ and $G_K$ are isomorphic.
\end{proposition}

\begin{proof}
Write $K=MHN^*$, where $M,N$ are monomial matrices with nonzero entries in $\mu_k$. Their block substitutions $\phi_k(M)$ and $\phi_k(N)$ are permutation matrices, and
\[
 \phi_k(K)=\phi_k(M)\phi_k(H)\phi_k(N)^{\top}.
\]
Consequently,
\[
 \mathcal M_K=(\phi_k(M)\oplus\phi_k(N))\mathcal M_H
 (\phi_k(M)\oplus\phi_k(N))^{\top}.
\]
This permutation similarity gives the required isomorphism.
\end{proof}

\begin{proposition}\label{prop:FAut}
Let $H\in\BH(n,k)$ and let $\ell$ be coprime to $k$. Then $G_H$ and $G_{H^{(\ell)}}$ are isomorphic.
\end{proposition}

\begin{proof}
Relabel every phase index by multiplication by $\ell$:
\[
 R_{i,a}\longmapsto R_{i,\ell a},\qquad
 C_{j,b}\longmapsto C_{j,\ell b}.
\]
If $h_{ij}=\zeta_k^{e_{ij}}$, then $b-a=e_{ij}$ is equivalent to
$\ell b-\ell a=\ell e_{ij}$ in $\mathbb Z_k$. The relabelling therefore preserves adjacency and loops, and gives an isomorphism from $G_H$ to $G_{H^{(\ell)}}$.
\end{proof}

Deleting the loops also removes the distinction between rows and columns. In fact,
\[
 R_{i,a}\longmapsto C_{i,-a},\qquad
 C_{j,b}\longmapsto R_{j,-b}
\]
is an isomorphism from the graph of $H$ to the graph of $H^{\top}$ with loops removed. With loops retained, an isomorphism must preserve the two parts.

\begin{proposition}\label{prop:Pform}
Let $H,K\in\BH(n,k)$. Every permutation matrix $\mathcal{P}$ satisfying
$\mathcal{P}\mathcal M_H\mathcal{P}^{\top}=\mathcal M_K$ has the form $\mathcal{P}=P\oplus Q$, where $P,Q$ are permutation matrices of order $nk$ satisfying
\[
 P\phi_k(H)Q^{\top}=\phi_k(K).
\]
\end{proposition}

\begin{proof}
The vertices in $\R$ are precisely the vertices carrying loops, so every isomorphism preserves $\R$ and $\C$. The stated identity follows by comparing the off-diagonal blocks of the adjacency matrices.
\end{proof}

\begin{theorem}[McKay \cite{McKay1979}]\label{thm:equivalence_iff_isomorphic}
Two real Hadamard matrices are equivalent if and only if their associated graphs are isomorphic.
\end{theorem}

\begin{remark}\label{remarkTranspose}
For real Hadamard matrices, the loopless graphs are isomorphic if and only if the matrices are equivalent after possibly transposing one of them \cite{McKay1979}. For general phase, Proposition~\ref{prop:FAut} shows that the graph also identifies field conjugates. Moreover, Proposition~\ref{prop:Pform} only gives arbitrary permutation matrices $P,Q$ of order $nk$; it does not require them to arise from monomial matrices of order $n$. Thus graph isomorphism alone does not establish monomial equivalence. This distinction is also relevant in algorithms for testing Butson equivalence; see \cite[Section~3]{BHPaper2}.
\end{remark}

\begin{example}\label{ex:isomorphic_inequivalent}
We give an example in exponent form. Let $H=(\zeta_{4}^{(E_H)_{ij}})$ and $K=(\zeta_{4}^{(E_K)_{ij}})$, where
\[
E_H=\begin{pmatrix}
0 & 0 & 0 & 0 & 0 & 0 & 0 & 0 & 0 & 0 & 0 & 0\\
0 & 0 & 0 & 0 & 0 & 0 & 2 & 2 & 2 & 2 & 2 & 2\\
0 & 0 & 1 & 2 & 2 & 3 & 0 & 0 & 1 & 2 & 2 & 3\\
0 & 0 & 1 & 2 & 2 & 3 & 2 & 2 & 3 & 0 & 0 & 1\\
0 & 1 & 3 & 0 & 2 & 2 & 0 & 1 & 3 & 0 & 2 & 2\\
0 & 1 & 3 & 0 & 2 & 2 & 2 & 3 & 1 & 2 & 0 & 0\\
0 & 2 & 1 & 3 & 0 & 2 & 0 & 2 & 1 & 3 & 0 & 2\\
0 & 2 & 1 & 3 & 0 & 2 & 2 & 0 & 3 & 1 & 2 & 0\\
0 & 2 & 3 & 2 & 1 & 0 & 0 & 2 & 3 & 2 & 1 & 0\\
0 & 2 & 3 & 2 & 1 & 0 & 2 & 0 & 1 & 0 & 3 & 2\\
0 & 3 & 2 & 1 & 3 & 1 & 0 & 3 & 2 & 1 & 3 & 1\\
0 & 3 & 2 & 1 & 3 & 1 & 2 & 1 & 0 & 3 & 1 & 3
\end{pmatrix},
\]
\[
E_K=\begin{pmatrix}
0 & 0 & 0 & 0 & 0 & 0 & 0 & 0 & 0 & 0 & 0 & 0\\
0 & 0 & 0 & 0 & 0 & 0 & 2 & 2 & 2 & 2 & 2 & 2\\
0 & 0 & 1 & 2 & 2 & 3 & 0 & 0 & 1 & 2 & 2 & 3\\
0 & 0 & 1 & 2 & 2 & 3 & 2 & 2 & 3 & 0 & 0 & 1\\
0 & 1 & 3 & 0 & 2 & 2 & 0 & 2 & 1 & 0 & 3 & 2\\
0 & 1 & 3 & 0 & 2 & 2 & 2 & 0 & 3 & 2 & 1 & 0\\
0 & 2 & 1 & 3 & 0 & 2 & 1 & 0 & 3 & 0 & 2 & 2\\
0 & 2 & 1 & 3 & 0 & 2 & 3 & 2 & 1 & 2 & 0 & 0\\
0 & 2 & 3 & 2 & 1 & 0 & 0 & 1 & 3 & 2 & 0 & 2\\
0 & 2 & 3 & 2 & 1 & 0 & 2 & 3 & 1 & 0 & 2 & 0\\
0 & 3 & 2 & 1 & 3 & 1 & 1 & 3 & 0 & 2 & 3 & 1\\
0 & 3 & 2 & 1 & 3 & 1 & 3 & 1 & 2 & 0 & 1 & 3
\end{pmatrix}.
\]
Both matrices belong to $\BH(12,4)$. Their graphs are isomorphic, although the matrices are inequivalent up to extended equivalence incorporating field automorphisms and transpose.
\end{example}

\subsection{Recovering a matrix from its graph}

The labelled fibres in Definition~\ref{def:BHGraph} record more information than an unlabelled graph. For prime phase, a common-neighbour condition characterises when a phase matrix recovered from these fibres is Hadamard. We require the following result of Lam and Leung \cite{LamLeung}.

\begin{lemma}\label{lem:prime_balance}
Let $k$ be prime and let $c_0,\ldots,c_{k-1}$ be nonnegative integers. Then
\[
 \sum_{a=0}^{k-1}c_a\zeta_k^a=0
 \quad\Longleftrightarrow\quad c_0=\cdots=c_{k-1}.
\]
\end{lemma}

\begin{proposition}\label{prop:reconstruction}
Let $k$ be prime. Let $G$ have labelled vertex set $\R\mathbin{\dot\cup}\C$ as in Definition~\ref{def:BHGraph}, with loops at all row vertices, no other edges within either part, and the following properties:
\begin{enumerate}[(i)]
\item each pair $R_i,C_j$ induces a perfect matching;
\item the simultaneous translation of all phase indices by $1$ preserves adjacency;
\item vertices in distinct row fibres have exactly $n/k$ common neighbours.
\end{enumerate}
Then $G=G_H$ for a unique $H\in\BH(n,k)$ in the given labelling. Conversely, every graph $G_H$ with $H\in\BH(n,k)$ has these properties.
\end{proposition}

\begin{proof}
By (i) and (ii), the matching between $R_i$ and $C_j$ is given by a unique shift $e_{ij}\in\mathbb Z_k$. Set $h_{ij}=\zeta_k^{e_{ij}}$. This recovers the adjacency rule \eqref{eq:adjacency_rule} and determines $H$ uniquely.

Fix distinct row indices $i,i'$. The number of common neighbours of $R_{i,a}$ and $R_{i',b}$ is
\begin{equation}\label{eq:common_neighbours}
 \#\{j:e_{ij}-e_{i'j}=b-a\}.
\end{equation}
Condition~(iii) says that every phase occurs $n/k$ times among the ratios $h_{ij}\overline{h_{i'j}}$. Their sum is therefore zero. Hence distinct rows are orthogonal and $HH^*=nI_n$.

Conversely, (i) and (ii) follow from the construction. If $H$ is Hadamard, orthogonality and Lemma~\ref{lem:prime_balance} show that every ratio occurs $n/k$ times. Equation~\eqref{eq:common_neighbours} then gives (iii).
\end{proof}

\section{Spectra of Butson Hadamard graphs}\label{section:spectra}

We now determine the spectrum of $G_H$. The key observation is that the regular representation separates into the characters of $\mu_k$, so the relevant matrices are the entrywise powers of $H$. Denote the $k \times k$ Fourier matrix by $\mathcal{F}_{k} = (\zeta_{k}^{ab})_{a,b \in \mathbb{Z}_{k}}$.

\begin{proposition}\label{prop:USpec}
For every $H\in\mu_k^{n\times n}$, the matrix $\phi_k(H)$ is unitarily similar to
\begin{equation}\label{eq:character_decomposition}
 \bigoplus_{\ell=0}^{k-1}H^{(\ell)}
 =J_n\oplus H\oplus\cdots\oplus H^{(k-1)}.
\end{equation}
In particular,
\[
 \spec(\phi_k(H))=\biguplus_{\ell=0}^{k-1}\spec(H^{(\ell)}).
\]
\end{proposition}

\begin{proof}
The unitary Fourier matrix $k^{-1/2}\mathcal{F}_k$ diagonalises $P_k$ \cite{Davis79}:
\[
 \frac{1}{k} \mathcal{F}_k^*P_k \mathcal{F}_k
 =\mathrm{diag}(1,\zeta_k,\ldots,\zeta_k^{k-1}).
\]
Conjugating $\phi_k(H)$ by $I_n\otimes(k^{-1/2} \mathcal{F}_k)$ replaces its $(i,j)$ block by $\diag(1,h_{ij},\ldots,h_{ij}^{k-1})$. Reordering the coordinates by character gives \eqref{eq:character_decomposition}.
\end{proof}

\begin{corollary}\label{cor:USpec}
Let $H,K\in\BH(n,k)$. If
\[
 \biguplus_{\ell=1}^{k-1}\spec(H^{(\ell)})
 =\biguplus_{\ell=1}^{k-1}\spec(K^{(\ell)}),
\]
then $\phi_k(H)$ and $\phi_k(K)$ are cospectral. If $k$ is prime, it is sufficient that $H$ and $K^{(r)}$ are cospectral for some $1\leq r<k$.
\end{corollary}

\begin{proof}
The first assertion follows from Proposition~\ref{prop:USpec}. For the second, equality of the spectra gives equality of the characteristic polynomials over $\mathbb Q(\zeta_k)$. Applying each field automorphism $\zeta_k\mapsto\zeta_k^\ell$ gives equality of the characteristic polynomials of $H^{(\ell)}$ and $K^{(r\ell)}$. Since multiplication by $r$ permutes the nonzero residues modulo the prime $k$, taking the multiset union proves the claim.
\end{proof}

The spectrum of $\mathcal M_H$ depends on the singular values of $\phi_k(H)$, rather than on its eigenvalues. The unitary form of Proposition~\ref{prop:USpec} makes this distinction explicit.

\begin{proposition}\label{prop:singular}
For $H\in\BH(n,k)$, the multiset of squared singular values of $\phi_k(H)$ is
\begin{equation}\label{eq:gram_spectra}
 \biguplus_{\ell=0}^{k-1}
 \spec\bigl(H^{(\ell)}(H^{(\ell)})^*\bigr).
\end{equation}
If $k$ is prime, these values are $n^2$ once, $0$ with multiplicity $n-1$, and $n$ with multiplicity $n(k-1)$.
\end{proposition}

\begin{proof}
Unitary similarity in Proposition~\ref{prop:USpec} preserves singular values, giving \eqref{eq:gram_spectra}. The block for $\ell=0$ is $J_nJ_n^{\top}=nJ_n$. If $k$ is prime, each nonzero $\ell$ gives a field automorphism, so $H^{(\ell)}(H^{(\ell)})^*=nI_n$.
\end{proof}

\begin{proposition}\label{prop:graph_spectrum}
Let $\sigma_1,\ldots,\sigma_{nk}$ be the singular values of $\phi_k(H)$, with multiplicities. Then the eigenvalues of $\mathcal M_H$ are
\begin{equation}\label{eq:graph_eigenvalues}
 \frac{1+\sqrt{1+4\sigma_j^2}}{2},\qquad
 \frac{1-\sqrt{1+4\sigma_j^2}}{2}
 \qquad(1\leq j\leq nk).
\end{equation}
The eigenvalues of the graph without loops are $\sigma_j$ and $-\sigma_j$ for the same indices.
\end{proposition}

\begin{proof}
Take a real singular value decomposition $U^{\top}\phi_k(H)V=\Sigma$. Conjugation by $U\oplus V$ changes \eqref{eq:AltAD} into
\[
 \begin{pmatrix}I_{nk}&\Sigma\\\Sigma&O_{nk}\end{pmatrix}.
\]
After reordering coordinates, this is the direct sum of the matrices
$\left(\begin{smallmatrix}1&\sigma_j\\\sigma_j&0\end{smallmatrix}\right)$, whose characteristic polynomials are $x^2-x-\sigma_j^2$. Deleting the loops replaces the diagonal entry $1$ by $0$ and gives eigenvalues $\pm\sigma_j$.
\end{proof}

\begin{theorem}\label{thm:prime-cospectral}
Let $k$ be prime and let $H\in\BH(n,k)$. The spectrum of $G_H$ is given in Table~\ref{EigenvalueMultiplicities}. In particular, all graphs associated with matrices in $\BH(n,k)$ are cospectral.
\end{theorem}

\begin{proof}
Combine Propositions~\ref{prop:singular} and \ref{prop:graph_spectrum}.
\end{proof}

\begin{table}[ht]\setlength\extrarowheight{6pt}
    \centering \begin{tabular}{ | l | l | l | l | l | l | l |}
    \hline
    Eigenvalue & $0$ & $1$ &  $\frac{1 + \sqrt{1 + 4n^2}}{2}$ &  $\frac{1 - \sqrt{1 + 4n^2}}{2}$ &  $\frac{1 + \sqrt{1 + 4n}}{2}$ &  $\frac{1 - \sqrt{1 + 4n}}{2}$ \\[2pt] \hline Multiplicity & $n-1$ & $n-1$ & $1$ & $1$ & $n(k-1)$ & $n(k-1)$ \\ \hline
    \end{tabular}
    \caption{Spectrum of $G_H$ for $H\in\BH(n,k)$ with $k$ prime.} \label{EigenvalueMultiplicities}
\end{table}

For composite $k$, the blocks with $\gcd(\ell,k)=1$ still equal $nI_n$ in \eqref{eq:gram_spectra}, but the remaining blocks need not do so. Thus $0$ and $1$ have multiplicity at least $n-1$, the two values $(1\pm\sqrt{1+4n^2})/2$ occur at least once each, and $(1\pm\sqrt{1+4n})/2$ each have multiplicity at least $n\varphi(k)$, where $\varphi$ is Euler's totient function. These contributions are added if values coincide. The nonunit exponents determine any additional contributions. Example~\ref{ex:k4} shows that the spectrum can then change under a switching that preserves the Butson property.

\section{WQH-switching and Orrick's Hadamard switchings}\label{section:new_results}

We now show that the two switching operations to construct Hadamard matrices in Section~\ref{subsection:orricks_switching} can be realised by WQH-switching on the associated Hadamard graphs. The switching sets are chosen from the phase copies of the fields in the relevant normal forms.

\subsection{Switching a closed quadruple}\label{subsection:ClosedQuadrupleWQH}

An arbitrary graph switching need not preserve the class of Hadamard graphs. The following example explains why a suitable switching partition is essential.

\begin{example}\label{ex:leaving_hadamard}
Let $H$ be the Sylvester Hadamard matrix of order $16$, indexed by $x,y\in\mathbb F_2^4$, with $h_{x,y}=(-1)^{x\cdot y}$. Write $e_1,e_2$ for the first two standard basis vectors. The sets
\[
 X^+=\{R_{0,0},R_{e_1,0}\},\qquad
 X^-=\{R_{e_2,1},R_{e_1+e_2,1}\}
\]
satisfy the one-pair WQH conditions: they contain only their own loops, and a column vertex has either one neighbour in each set or two in one set and none in the other.

The all-or-none case occurs precisely for the four columns with $y_1=y_2=0$. Before switching, $R_{0,0}$ and $R_{e_2,0}$ have eight common neighbours. The switching removes four of these and adds no new common neighbours to this pair. They therefore have four common neighbours in the switched graph. In a Hadamard graph of order $16$, distinct looped vertices have either zero common neighbours, if they belong to the same row fibre, or eight, if they belong to different fibres. Hence the switched graph is not isomorphic to any Hadamard graph of order $16$.
\end{example}

\begin{theorem}\label{thm:ClosedQuadrupleWQH}
Let $H$ be a Hadamard matrix of order $n>4$ with a closed quadruple, and let $H'$ be obtained by switching that quadruple. Then $G_{H'}$ is isomorphic to a graph obtained from $G_H$ by WQH-switching.
\end{theorem}

\begin{proof}
By Theorems~\ref{thm:orrick_quad} and \ref{thm:equivalence_iff_isomorphic}, it suffices to take $H$ in the form \eqref{closedquad} and negate its initial $4\times m$ block, where $m=n/4$. Let
\[
 X^+=\{C_{j,0}:0\leq j<m\},\qquad
 X^-=\{C_{j,1}:0\leq j<m\}.
\]
These are equally sized independent sets, with no edges between them. For $a\in\mathbb Z_2$, put $U^a=\{R_{i,a}:0\leq i<4\}$. Every vertex of $U^0$ is adjacent to all of $X^+$ and none of $X^-$; every vertex of $U^1$ has the reverse neighbourhood.

Every row of $A_1$ in \eqref{closedquad} has sum zero, so every other row vertex has $m/2$ neighbours in each of $X^+$ and $X^-$. The remaining column vertices have no neighbours in either set. The one-pair WQH conditions therefore hold. Switching interchanges exactly the matchings between the first four row fibres and the first $m$ column fibres. By \eqref{eq:adjacency_rule}, this negates the initial $4\times m$ block and gives $G_{H'}$.
\end{proof}

\subsection{Switching a Hall set}\label{subsection:HallSetWQH}

\begin{theorem}\label{thm:HallSetWQH}
Let $H$ be a Hadamard matrix of order $n\geq12$, with $n\equiv4\pmod8$, containing a Hall set. If $H'$ is obtained by switching that Hall set, then $G_{H'}$ is isomorphic to a graph obtained from $G_H$ by WQH-switching.
\end{theorem}

\begin{proof}
Use the normal form \eqref{Hallsetmatrixblocks}, with $m=(n-4)/4$. By Theorems~\ref{thm:orrick_hall} and \ref{thm:equivalence_iff_isomorphic}, it suffices to negate $F_1$ and $G_1$. Define
\begin{align*}
 X_{1}^+&=\{C_{j,0}:4\leq j<m+4\},&
 X_{1}^-&=\{C_{j,1}:4\leq j<m+4\},\\
 X_{2}^+&=\{R_{i,0}:4\leq i<m+4\},&
 X_{2}^-&=\{R_{i,1}:4\leq i<m+4\},
\end{align*}
and let $D$ consist of the remaining vertices. The four sets above have size $m$.

We check condition~(ii) of Theorem~\ref{thm:WQH} for the two pairs $(X_{1}^+,X_{1}^-)$ and $(X_{2}^+,X_{2}^-)$. There are no edges within the column part, while every row vertex has its own loop and no other row neighbour. The corresponding difference constants are therefore $c_{11}=0$ and $c_{22}=1$. Since every row and column of $A_{11}$ has sum $2$, each contains $m/2+1$ positive entries and $m/2-1$ negative entries. Hence the two cross-difference constants are $c_{12}=c_{21}=2$.

For $a\in\mathbb Z_2$, put
\[
 U^a=\{R_{i,a}:0\leq i<4\},\qquad
 V^a=\{C_{j,a}:0\leq j<4\}.
\]
The all-one blocks $F_1,G_1$ imply that $U^0$ is complete to $X_{1}^+$ and anticomplete to $X_{1}^-$, while $U^1$ has the reverse neighbourhood. Similarly, $V^0$ is complete to $X_{2}^+$ and anticomplete to $X_{2}^-$, and $V^1$ has the reverse neighbourhood. These vertices have no neighbours in the other pair of switching sets.

All other vertices of $D$ are balanced on each pair. Indeed, the relevant off-diagonal blocks $A_{a1}$ and $A_{1a}$ have row and column sums zero, giving $m/2$ neighbours in each set of the appropriate pair; the counts for the other pair are both zero. Thus all WQH conditions hold.

The switching changes exactly the adjacencies between $U^0\cup U^1$ and $X_{1}^+\cup X_{1}^-$, and between $V^0\cup V^1$ and $X_{2}^+\cup X_{2}^-$. These changes negate $F_1$ and $G_1$, so the resulting graph is $G_{H'}$.
\end{proof}

\section{Switching methods for Butson matrices and their graphs}\label{section:butson_switching}

Theorem~\ref{thm:prime-cospectral} shows that a switching which stays within $\BH(n,k)$ automatically preserves the associated graph spectrum when $k$ is prime. For composite phase, this conclusion requires further hypotheses. In this section we first recall two matrix switchings from \cite{GenSwitching}, then characterise the GM and WQH reflections on Butson matrices and give graph realisations of these operations and of balanced block switching.

\subsection{Block switching and composite-phase examples}

The following form of block switching uses separate parameters for the number of selected rows and the phase. Its short proof also explains the balance condition needed later.

\begin{proposition}[Block switching; cf.\ {\cite[Proposition~5.1]{GenSwitching}}]\label{prop:FourierSwitch}
Let $q,m\geq1$, put $n=qm$, and let $H\in\BH(n,k)$ have its first $q$ rows equal to
\[
 K_0\otimes\one_m^{\top},\qquad K_0\in\BH(q,k).
\]
Multiplying any one of the $q\times m$ blocks in these rows by an element of $\mu_k$ produces another matrix in $\BH(n,k)$.
\end{proposition}

\begin{example}\label{ex:BH66}
Let $\omega=\zeta_3$. Consider the following matrices in $\BH(6,6)$:
\[
H_1=\begin{pmatrix}
1 & 1 & 1 & 1 & 1 & 1\\
1 & 1 & \omega & \omega & \omega^2 & \omega^2\\
1 & 1 & \omega^2 & \omega^2 & \omega & \omega\\
1 & -1 & -1 & 1 & 1 & -1\\
1 & -1 & -\omega & \omega & \omega^2 & -\omega^2\\
1 & -1 & -\omega^2 & \omega^2 & \omega & -\omega
\end{pmatrix},\qquad
H_2=\begin{pmatrix}
1 & 1 & \omega & \omega & 1 & 1\\
1 & 1 & \omega^2 & \omega^2 & \omega^2 & \omega^2\\
1 & 1 & 1 & 1 & \omega & \omega\\
1 & -1 & -1 & 1 & 1 & -1\\
1 & -1 & -\omega & \omega & \omega^2 & -\omega^2\\
1 & -1 & -\omega^2 & \omega^2 & \omega & -\omega
\end{pmatrix}.
\]
The first three rows of $H_1$ are $\mathcal{F}_3\otimes\one_2^{\top}$, and $H_2$ is obtained by multiplying their middle $3\times2$ block by $\omega$. Proposition~\ref{prop:FourierSwitch} therefore applies. We can computationally check that $G_{H_{1}}$ and $G_{H_{2}}$ are non-isomorphic cospectral mates.
\end{example}

The next result extends the block pattern underlying Hall-set switching.

\begin{theorem}[{\cite[Theorem~5.3 and Corollary~5.4]{GenSwitching}}]\label{thm:genHall}
Let $m\geq1$, let $n=k(m+1)$, and let $H\in\BH(n,k)$ have the form
\[
 H=\begin{pmatrix}
 S& \mathcal{F}_k\otimes\one_m^{\top}\\
 \mathcal{F}_k^*\otimes\one_m&A
 \end{pmatrix},
\]
where $S\in\BH(k,k)$ is circulant and $A$ has order $km$. Choose an index $a\in\{0,\ldots,k-1\}$. Multiply the $k\times m$ block in the first $k$ rows and columns $k+am,\ldots,k+(a+1)m-1$ by $\alpha\in\mu_k$. Simultaneously multiply the $m\times k$ block in the corresponding rows and first $k$ columns by $\overline\alpha$. The resulting matrix belongs to $\BH(n,k)$.
\end{theorem}

For composite phase, this operation can change the spectrum of the associated graph.

\begin{example}\label{ex:k4}
We use the matrix from \cite[Example~5.8]{GenSwitching}, $H = (\zeta_{4}^{e_{ij}})$, where
\[
E=\left(\begin{array}{cccc|cc|cc|cc|cc}
0 & 1 & 2 & 1 & 0 & 0 & 0 & 0 & 0 & 0 & 0 & 0 \\
1 & 0 & 1 & 2 & 0 & 0 & 1 & 1 & 2 & 2 & 3 & 3 \\
2 & 1 & 0 & 1 & 0 & 0 & 2 & 2 & 0 & 0 & 2 & 2 \\
1 & 2 & 1 & 0 & 0 & 0 & 3 & 3 & 2 & 2 & 1 & 1 \\
\hline
0 & 0 & 0 & 0 & 1 & 1 & 0 & 2 & 0 & 2 & 0 & 2 \\
0 & 0 & 0 & 0 & 1 & 1 & 2 & 0 & 2 & 0 & 2 & 0 \\
\hline
0 & 3 & 2 & 1 & 0 & 2 & 2 & 2 & 0 & 2 & 2 & 0 \\
0 & 3 & 2 & 1 & 2 & 0 & 2 & 2 & 2 & 0 & 0 & 2 \\
\hline
0 & 2 & 0 & 2 & 0 & 2 & 2 & 0 & 3 & 3 & 0 & 2 \\
0 & 2 & 0 & 2 & 2 & 0 & 0 & 2 & 3 & 3 & 2 & 0 \\
\hline
0 & 1 & 2 & 3 & 0 & 2 & 0 & 2 & 2 & 0 & 2 & 2 \\
0 & 1 & 2 & 3 & 2 & 0 & 2 & 0 & 0 & 2 & 2 & 2
\end{array}\right).
\]
It belongs to $\BH(12,4)$. Let $K$ be obtained by multiplying the initial all-one $4\times2$ block by $\zeta_{4}$ and the corresponding $2\times4$ block by $-\zeta_{4}$. By Theorem~\ref{thm:genHall}, $K\in\BH(12,4)$. Computationally we verify that $G_{H}$ and $G_{K}$ are not cospectral.
\end{example}

\subsection{The GM reflection on Butson matrices}\label{subsection:butson_gm_reflection}

We now use the GM reflection $Q_m$ from \eqref{eq:butson_gm_Q} as a matrix-switching operation. Orthogonality is preserved automatically. The first question is therefore when the transformed entries remain in $\mu_k$. Related extensions of GM-switching to complex unit gain graphs appear in \cite{BBCD2021paper2}; here we work with Butson matrices and their associated unweighted graphs.

For $x\in\mathbb C^m$,
\begin{equation}\label{eq:butson_gm_action}
 Q_mx=\frac{2}{m}\left(\sum_{r=0}^{m-1}x_r\right)\one_m-x.
\end{equation}

\begin{lemma}\label{lem:butson_gm_columns}
Let $x\in\mu_k^m$. Then $Q_mx\in\mu_k^m$ if and only if at least one of the following conditions holds:
\begin{enumerate}[(a)]
\item all entries of $x$ are equal;
\item $m$ is even, and $x$ contains two distinct roots, each occurring $m/2$ times;
\item $k$ is even and $\sum_r x_r=0$.
\end{enumerate}
In these cases, $Q_m$ fixes $x$, interchanges its two values, or negates $x$, respectively. The descriptions agree when the conditions overlap.
\end{lemma}

\begin{proof}
Sufficiency follows directly from \eqref{eq:butson_gm_action}. For necessity, put $z=(2/m)\sum_r x_r$. If $z=0$, then $Q_mx=-x$, and $-x_r\in\mu_k$ holds precisely when $-1\in\mu_k$, or equivalently when $k$ is even.

Suppose $z\neq0$. Every entry $\alpha$ of $x$ must satisfy $|\alpha|=|z-\alpha|=1$. The two circles with centres $0$ and $z$ and radius $1$ have at most two intersection points. Thus $x$ contains at most two distinct values. If it is constant, (a) holds. Otherwise, let its values be $\alpha\neq\beta$, with $\alpha$ occurring a proportion $p$ of the time. Then $z=2(p\alpha+(1-p)\beta)$. The identity $|z-\alpha|^2=|z-\beta|^2$ gives
\[
 (2p-1)|\alpha-\beta|^2=0.
\]
Hence $p=1/2$, proving (b).
\end{proof}

\begin{theorem}\label{thm:butson_gm_matrix}
Let $H\in\BH(n,k)$, and let $S_1,\ldots,S_t$ be pairwise disjoint nonempty sets of row indices. Write $s_\ell=|S_\ell|$, and let $U$ act as $Q_{s_\ell}$ on the coordinates indexed by $S_\ell$ and as the identity on the remaining coordinates. Set $K=UH$, so
\begin{equation}\label{eq:butson_gm_switch}
 k_{ij}=\begin{cases}
 \displaystyle\frac{2}{s_\ell}\sum_{r\in S_\ell}h_{rj}-h_{ij},&i\in S_\ell,\\[5pt]
 h_{ij},&i\notin\bigcup_{\ell=1}^t S_\ell.
 \end{cases}
\end{equation}
Then $K\in\BH(n,k)$ if and only if every restricted column $(h_{ij})_{i\in S_\ell}$ satisfies a condition of Lemma~\ref{lem:butson_gm_columns}. Whenever this holds, the operation is an involution on the same row sets.
\end{theorem}

\begin{proof}
After reordering the row coordinates, $U$ is a direct sum of the matrices $Q_{s_\ell}$ and an identity matrix. Hence $U=U^*$ and $U^2=I_n$. Therefore
\[
 KK^*=UHH^*U^*=nI_n.
\]
The remaining requirement is that all entries of $K$ lie in $\mu_k$, which is exactly the condition of Lemma~\ref{lem:butson_gm_columns}. Finally, $UK=H$.
\end{proof}

\begin{remark}\label{remark:phased_GM}
The same construction can be applied after multiplying the rows by prescribed phases. If $T=\diag(d_0,\ldots,d_{n-1})$, with $d_i\in\mu_k$, satisfies the required entry conditions for $TH$, then undoing the phases gives
\begin{equation}\label{eq:butson_gm_phased}
 K=T^*UTH\in\BH(n,k).
\end{equation}
Column versions follow by transposition.
\end{remark}

The constant and equally split two-value cases have a direct GM interpretation on the associated graph. The zero-sum condition alone does not give the same graph partition.

\begin{theorem}\label{thm:butson_gm_lift}
Let $H\in\BH(n,k)$, and let $S_1,\ldots,S_t$ be disjoint row sets of even sizes $s_\ell\geq2$. Suppose that every column restricted to $S_\ell$ is constant or contains two distinct roots, each occurring $s_\ell/2$ times. Define $U$ as in Theorem~\ref{thm:butson_gm_matrix} and set $K=UH$. Then $K\in\BH(n,k)$ and $G_K$ is obtained from $G_H$ by GM-switching with cells
\begin{equation}\label{eq:butson_gm_cells}
 X_{\ell,a}=\{R_{i,a}:i\in S_\ell\},\qquad
 1\leq\ell\leq t,\quad a\in\mathbb Z_k,
\end{equation}
and remainder $D=V(G_H)\setminus\bigcup_{\ell,a}X_{\ell,a}$.
Conversely, the GM conditions for these cells force each restricted column to have one of the two stated forms.
\end{theorem}

\begin{proof}
The matrix assertion follows from Theorem~\ref{thm:butson_gm_matrix}. The cells consist of row vertices, so every vertex has just its own loop within its cell and no neighbours in any other cell. This gives condition~(i) of Theorem~\ref{thm:butson_gm_classical}. Remaining row vertices have no neighbours in the switching cells. For a column vertex,
\begin{equation}\label{eq:butson_gm_count}
 d_{X_{\ell,a}}(C_{j,b})
 =\#\{i\in S_\ell:h_{ij}=\zeta_k^{b-a}\}.
\end{equation}
The column hypothesis makes this number $0$, $s_\ell/2$, or $s_\ell$, giving condition~(ii).

A constant restricted column is unchanged. If two roots occur equally often, the switching complements the two half-neighbourhoods and interchanges their positions. This is exactly the action of $Q_{s_\ell}$ in Lemma~\ref{lem:butson_gm_columns}, so the resulting graph is $G_K$.

Conversely, fix $\ell,j,b$ and vary $a$ in \eqref{eq:butson_gm_count}. The counts are the phase multiplicities in the restricted column and sum to $s_\ell$. If every count lies in $\{0,s_\ell/2,s_\ell\}$, either one phase occurs $s_\ell$ times or two phases occur $s_\ell/2$ times each.
\end{proof}

In the fixed vertex order \eqref{eq:partition}, the conjugating matrix is particularly simple. Put $\widehat U=U\otimes I_k$ and $\mathcal Q=\widehat U\oplus I_{nk}$. Then
\begin{equation}\label{eq:butson_gm_liftidentity}
 \phi_k(K)=\widehat U\phi_k(H),\qquad
 \mathcal M_K=\mathcal Q\mathcal M_H\mathcal Q^{\top}.
\end{equation}
These identities follow blockwise from the constant or equally split two-value condition. In particular, $G_H$ and $G_K$ are cospectral for every phase $k$ covered by the theorem, including odd and composite phases.

\begin{example}\label{ex:butson_gm_ternary}
Let $\omega=\zeta_3$, and consider
\[
H=\begin{pmatrix}
1 & 1 & 1 & 1 & 1 & 1\\
1 & 1 & \omega & \omega & \omega^2 & \omega^2\\
\omega & \omega^2 & \omega & 1 & 1 & \omega^2\\
\omega & \omega^2 & 1 & \omega & \omega^2 & 1\\
1 & \omega^2 & \omega^2 & \omega & 1 & \omega\\
1 & \omega^2 & \omega & \omega^2 & \omega & 1
\end{pmatrix}.
\]
Then $H\in\BH(6,3)$. On the first four rows, each restricted column contains two distinct phases twice each. Taking $U=Q_4\oplus I_2$, we obtain
\[
K=UH=\begin{pmatrix}
\omega & \omega^2 & \omega & \omega & \omega^2 & \omega^2\\
\omega & \omega^2 & 1 & 1 & 1 & 1\\
1 & 1 & 1 & \omega & \omega^2 & 1\\
1 & 1 & \omega & 1 & 1 & \omega^2\\
1 & \omega^2 & \omega^2 & \omega & 1 & \omega\\
1 & \omega^2 & \omega & \omega^2 & \omega & 1
\end{pmatrix}\in\BH(6,3).
\]
By Theorem~\ref{thm:butson_gm_lift}, the graph $G_K$ is obtained from $G_H$ by GM-switching on the three cells
\[
 X_a=\{R_{0,a},R_{1,a},R_{2,a},R_{3,a}\},\qquad a\in\mathbb Z_3.
\]
This operation interchanges phases within each restricted column and applies at odd phase. This does not assert that $H$ and $K$ are inequivalent.
\end{example}

\subsection{The WQH reflection on Butson matrices}\label{subsection:butson_wqh_reflection}

We next apply the paired reflection underlying WQH-switching to Butson matrices. Extensions of WQH-switching to complex unit and general gain graphs appear in \cite{ABK2024,BBCD2021paper2}. We distinguish the condition that this reflection preserves the allowed matrix entries from the stronger phase-multiplicity condition that gives WQH-switching on the associated unweighted graphs.

For $s\geq1$, put $z=(\one_s^{\top},-\one_s^{\top})^{\top}$ and define
\begin{equation}\label{eq:butson_wqh_W}
 W_s=I_{2s}-\frac1szz^{\top}
 =\begin{pmatrix}
 I_s-\frac1sJ_s&\frac1sJ_s\\
 \frac1sJ_s&I_s-\frac1sJ_s
 \end{pmatrix}.
\end{equation}
This is the one-pair WQH matrix from Section~\ref{subsection:graph_switching}. Since $z^{\top}z=2s$ and $z^{\top}\one_{2s}=0$, we have $W_s^{\top}=W_s$, $W_s^2=I_{2s}$ and $W_s\one_{2s}=\one_{2s}$. For $u,v\in\mathbb C^s$, its action is
\begin{equation}\label{eq:butson_wqh_action}
 W_s\binom{u}{v}=\binom{u-\delta\one_s}{v+\delta\one_s},
 \qquad
 \delta=\frac1s\left(\sum_{r=0}^{s-1}u_r-\sum_{r=0}^{s-1}v_r\right).
\end{equation}

\begin{lemma}\label{lem:butson_wqh_columns}
Let $u,v\in\mu_k^s$. Then $W_s(u^{\top},v^{\top})^{\top}\in\mu_k^{2s}$ if and only if at least one of the following conditions holds:
\begin{enumerate}[(a)]
\item $\sum_r u_r=\sum_r v_r$;
\item both $u$ and $v$ are constant, possibly with different values;
\item $k$ is even, and the signed vector $(u^{\top},-v^{\top})^{\top}$ contains two distinct roots, each occurring $s$ times.
\end{enumerate}
In case~(a), the reflection fixes the column. In case~(b), it interchanges the two constant values. In case~(c), if the two signed values are $\alpha,\beta$, then $\delta=\alpha+\beta$ and the signed values are replaced by $-\beta,-\alpha$, respectively. These descriptions agree when the conditions overlap. In particular, for odd $k$, conditions~(a) and~(b) are the only possibilities.
\end{lemma}

\begin{proof}
Sufficiency follows from \eqref{eq:butson_wqh_action}. For necessity, the case $\delta=0$ is (a). Suppose $\delta\neq0$, and write $y=(u^{\top},-v^{\top})^{\top}$. Every entry $\gamma$ of $y$ satisfies $|\gamma|=|\delta-\gamma|=1$. The two circles have distinct centres, so $y$ has at most two distinct values. If $y$ is constant, (b) holds. Otherwise, let its two values be $\alpha\neq\beta$, with $\alpha$ occurring a proportion $p$ of the time. Since
\[
 \delta=\frac1s\sum_r y_r=2\bigl(p\alpha+(1-p)\beta\bigr),
\]
the equality $|\delta-\alpha|^2=|\delta-\beta|^2$ gives
\[
 (2p-1)|\alpha-\beta|^2=0.
\]
Thus each value occurs $s$ times. If $k$ is even, all entries of $y$ lie in $\mu_k$, giving (c). If $k$ is odd, the sets $\mu_k$ and $-\mu_k$ are disjoint. The entries of $u$ lie in the first set and those of $-v$ in the second. With only two values available, each half must therefore be constant, giving (b).
\end{proof}

\begin{theorem}\label{thm:butson_wqh_matrix}
Let $H\in\BH(n,k)$, and let $S_1^+,S_1^-,\ldots,S_t^+,S_t^-$ be pairwise disjoint nonempty row sets with $|S_\ell^+|=|S_\ell^-|=s_\ell$. In $\mathbb R^n$, put
\[
 z_\ell=\one_{S_\ell^+}-\one_{S_\ell^-},\qquad
 V=I_n-\sum_{\ell=1}^t\frac1{s_\ell}z_\ell z_\ell^{\top},
 \qquad K=VH
\]
where $z_{\ell}$ is the difference of characteristic vectors. Equivalently, with
\[
 \delta_{\ell j}=\frac1{s_\ell}
 \left(\sum_{r\in S_\ell^+}h_{rj}-\sum_{r\in S_\ell^-}h_{rj}\right),
\]
the entries of $K$ are
\begin{equation}\label{eq:butson_wqh_switch}
 k_{ij}=\begin{cases}
 h_{ij}-\delta_{\ell j},&i\in S_\ell^+,\\
 h_{ij}+\delta_{\ell j},&i\in S_\ell^-,\\
 h_{ij},&i\notin\bigcup_{\ell=1}^t(S_\ell^+\cup S_\ell^-).
 \end{cases}
\end{equation}
Then $K\in\BH(n,k)$ if and only if every column restricted to each ordered pair $(S_\ell^+,S_\ell^-)$ satisfies Lemma~\ref{lem:butson_wqh_columns}. Whenever this holds, the operation is an involution on the same row sets.
\end{theorem}

\begin{proof}
The vectors $z_\ell$ have disjoint supports and squared norms $2s_\ell$. Hence $V^{\top}=V$ and $V^2=I_n$, and
\[
 KK^*=VHH^*V^{\top}=nI_n.
\]
The remaining requirement is entry preservation, which is precisely Lemma~\ref{lem:butson_wqh_columns} on each pair. The identity $VK=H$ proves involutivity.
\end{proof}

\begin{remark}\label{remark:butson_wqh_gm_relation}
The relation with the preceding construction is
\begin{equation}\label{eq:butson_wqh_signed_gm}
 W_s=-D_sQ_{2s}D_s,\qquad D_s=\diag(I_s,-I_s).
\end{equation}
For even $k$, these row signs are allowed phase scalings, so the WQH reflection is related to the phased GM reflection of Remark~\ref{remark:phased_GM} by such scalings. The paired graph partition below gives its direct WQH interpretation. For odd $k$, the signs in $D_s$ do not preserve the specified phase group. The entry criterion may also be applied after prescribed row phase scalings, as in Remark~\ref{remark:phased_GM}, to obtain phased versions. Transposition gives column versions.
\end{remark}

Equality of the two column sums suffices to fix a column in the matrix operation. To fix all corresponding graph incidences, we require equality of the individual phase multiplicities.

\begin{theorem}\label{thm:butson_wqh_lift}
Let $H$, the paired row sets and $V$ be as in Theorem~\ref{thm:butson_wqh_matrix}. Suppose that for every pair and every column either
\begin{enumerate}[(a)]
\item each element of $\mu_k$ occurs equally often on $S_\ell^+$ and $S_\ell^-$; or
\item the column is constant on each of $S_\ell^+$ and $S_\ell^-$, possibly with different values.
\end{enumerate}
Then $K=VH\in\BH(n,k)$, and $G_K$ is obtained from $G_H$ by WQH-switching with paired cells
\begin{equation}\label{eq:butson_wqh_cells}
 X_{\ell,a}^{\pm}=\{R_{i,a}:i\in S_\ell^{\pm}\},
 \qquad 1\leq\ell\leq t,\quad a\in\mathbb Z_k,
\end{equation}
and remainder $D$ consisting of all other vertices. Conversely, the WQH conditions for these paired cells force each restricted column to have one of the two stated forms.
\end{theorem}

\begin{proof}
The matrix assertion follows from Lemma~\ref{lem:butson_wqh_columns} and Theorem~\ref{thm:butson_wqh_matrix}. Each cell has size $s_\ell$. Within their union, every vertex has only its own loop, so the difference constants in Theorem~\ref{thm:WQH}(ii) are $1$ for a pair with itself and $0$ for distinct pairs. Remaining row vertices have no neighbours in the cells. For a column vertex, the neighbourhood counts are
\begin{equation}\label{eq:butson_wqh_count}
 d_{X_{\ell,a}^{\pm}}(C_{j,b})
 =\#\{i\in S_\ell^{\pm}:h_{ij}=\zeta_k^{b-a}\}.
\end{equation}
In case~(a), these counts agree for the two cells, and all incidences are fixed. In case~(b), equal constant values also give equal counts. Distinct constant values give, for the two relevant phases, all of one cell and none of the other. WQH-switching interchanges these neighbourhoods, exactly interchanging the two constant values in \eqref{eq:butson_wqh_action}. Thus the switched graph is $G_K$.

Conversely, fix $\ell,j$, and let $m_r^{\pm}$ be the multiplicity of $\zeta_k^r$ on $S_\ell^{\pm}$. Varying $a$ in \eqref{eq:butson_wqh_count}, the WQH condition requires, for every $r$, either $m_r^+=m_r^-$ or
\[
 (m_r^+,m_r^-)\in\{(s_\ell,0),(0,s_\ell)\}.
\]
If any pair of counts differs, one half is constant. Since the multiplicities in each half sum to $s_\ell$, applying the same condition to all other phases forces the other half to be constant too. Otherwise all phase multiplicities agree.
\end{proof}

In the fixed vertex order \eqref{eq:partition}, put $\widehat V=V\otimes I_k$ and $\mathcal W=\widehat V\oplus I_{nk}$. Under the hypotheses of Theorem~\ref{thm:butson_wqh_lift}, we have
\begin{equation}\label{eq:butson_wqh_liftidentity}
 \phi_k(K)=\widehat V\phi_k(H),\qquad
 \mathcal M_K=\mathcal W\mathcal M_H\mathcal W^{\top}.
\end{equation}
Indeed, equal phase multiplicities give equal sums of the permutation blocks on the two halves, while constant halves have their permutation blocks interchanged. This verifies the first identity directly, without treating $\phi_k$ as additive on complex numbers. The second identity proves cospectrality for every phase covered by the lifting criterion.

\begin{corollary}\label{cor:butson_wqh_prime}
Let $k$ be prime, and let $H,V$ be as in Theorem~\ref{thm:butson_wqh_matrix}. If $VH\in\BH(n,k)$, then the hypotheses of Theorem~\ref{thm:butson_wqh_lift} hold. Thus every admissible paired reflection at prime phase has the stated WQH realisation on the associated graph.
\end{corollary}

\begin{proof}
For odd $k$, Lemma~\ref{lem:butson_wqh_columns} allows only equal sums or constant halves. For $k=2$, its case~(c) also has equal sums, since the two distinct signed values are $1$ and $-1$. It therefore remains to consider equal sums. Write $m_r^{\pm}$ for the phase multiplicities of the two halves and put
\[
 f(x)=\sum_{r=0}^{k-1}(m_r^+-m_r^-)x^r.
\]
We have $f(\zeta_k)=0$. Since $k$ is prime and $\deg f\leq k-1$, the minimal polynomial of $\zeta_k$ gives $f=c(1+x+\cdots+x^{k-1})$ for some rational $c$. The two halves have equal size, so $f(1)=0$, whence $c=0$. All phase multiplicities therefore agree.
\end{proof}

The corollary identifies a specific graph switching, in addition to the cospectrality already guaranteed by Theorem~\ref{thm:prime-cospectral}. For composite phase, the distinction between sums and phase multiplicities must be retained. For example, in $H=(\zeta_{4}^{ij})_{0\leq i,j<4}\in\BH(4,4)$, take $S^+=\{0,2\}$ and $S^-=\{1,3\}$. The reflection fixes columns $0,1,3$ and negates column $2$, so it preserves the Butson property. However, in column $1$ the two halves are $(1,-1)^{\top}$ and $(\zeta_{4},-\zeta_{4})^{\top}$. Their sums agree but their phase multiplicities do not, and $C_{1,0}$ has respectively $1$ and $0$ neighbours in $X_{1,0}^{+}$ and $X_{1,0}^{-}$. These cells fail the WQH conditions. In this example the matrices are column-equivalent; the obstruction concerns this particular graph partition.

\begin{example}\label{ex:butson_wqh_ternary}
Let $H\in\BH(6,3)$ be the matrix in Example~\ref{ex:butson_gm_ternary}, and take $S^+=\{0,1\}$ and $S^-=\{2,3\}$. Its first two columns are constant on each half, with values $(1,\omega)$ and $(1,\omega^2)$, respectively. In each of the other four columns, the two halves have identical phase multisets. Thus $V=W_2\oplus I_2$ gives
\[
K=VH=\begin{pmatrix}
\omega & \omega^2 & 1 & 1 & 1 & 1\\
\omega & \omega^2 & \omega & \omega & \omega^2 & \omega^2\\
1 & 1 & \omega & 1 & 1 & \omega^2\\
1 & 1 & 1 & \omega & \omega^2 & 1\\
1 & \omega^2 & \omega^2 & \omega & 1 & \omega\\
1 & \omega^2 & \omega & \omega^2 & \omega & 1
\end{pmatrix}\in\BH(6,3).
\]
Only the first two columns on the four selected rows change. Theorem~\ref{thm:butson_wqh_lift} realises this operation as WQH-switching on the three pairs
\[
 X_a^+=\{R_{0,a},R_{1,a}\},\qquad
 X_a^-=\{R_{2,a},R_{3,a}\},\qquad a\in\mathbb Z_3.
\]
As in Example~\ref{ex:butson_gm_ternary}, no inequivalence or graph nonisomorphism is asserted.
\end{example}

\subsection{Cyclic graph switching}\label{subsection:cyclic_switching}

For real matrices, negating a switched block interchanges the two phase classes of the associated graph. Multiplying a Butson block by $\zeta_k$ instead rotates $k$ phase classes. We make this observation precise in a graph-switching construction.

\begin{theorem}\label{thm:cyclic_switching}
Let $G$ be a graph with disjoint vertex sets $X_0,\ldots,X_{k-1}$ of the same size $m\geq1$, where $k\geq2$. Suppose that their union $X$ is independent and has no loops. For every vertex $v\notin X$, suppose that either
\begin{enumerate}[(i)]
\item $N(v)\cap X=X_a$ for some $a\in\mathbb Z_k$; or
\item $d_{X_0}(v)=\cdots=d_{X_{k-1}}(v)$.
\end{enumerate}
For each vertex in case~(i), replace its neighbourhood $X_a$ by $X_{a+1}$, with indices modulo $k$. Leave all other adjacencies unchanged. The resulting graph is cospectral with $G$.
\end{theorem}

\begin{proof}
Order the vertices of $X$ by the classes $X_0,\ldots,X_{k-1}$, and put
\[
 E=\frac1mJ_m,\qquad
 Q=I_k\otimes(I_m-E)+P_k\otimes E.
\]
Since $E$ and $I_m-E$ are orthogonal projections and $P_k$ is orthogonal, $QQ^{\top}=I_{km}$. Let $b$ be the row of incidences from a vertex outside $X$ to $X$. If its $k$ blocks have equal sums, then $bQ=b$. If $b$ is all ones in block $a$ and zero elsewhere, multiplication by $Q$ moves that all-one block to $a+1$.

Write the adjacency matrix with the vertices outside $X$ first. Conjugating by $I_{|V(G)\setminus X|}\oplus Q$ preserves the zero block on $X$, preserves the adjacency block outside $X$, and changes the off-diagonal blocks exactly as stated. It therefore produces the adjacency matrix of the switched graph. Orthogonal similarity proves cospectrality.
\end{proof}

When $k=2$, this is the one-pair WQH operation for independent switching sets. No divisibility condition on $m$ is required by the graph theorem.

We next apply the theorem to Proposition~\ref{prop:FourierSwitch}. After scaling the first $q$ rows, we may assume that the selected $q\times m$ block is $J_{q,m}$. Let $A$ be the submatrix formed by the remaining rows and those same $m$ columns. The proof of Proposition~\ref{prop:FourierSwitch} gives zero row sums for $A$. To obtain the graph switching, we require the stronger condition: that each element of $\mu_k$ occurs equally often in every row of $A$.
If $A$ has no rows, this condition is vacuous. Otherwise it implies $k\mid m$. For prime $k$, it follows from the zero row sums by Lemma~\ref{lem:prime_balance}.

\begin{corollary}\label{cor:cyclic_BH}
Under the hypotheses of Proposition~\ref{prop:FourierSwitch}, assume that each element of $\mu_k$ occurs equally often in the submatrix $A$ completing the columns of the selected block, normalised to $J_{q,m}$. Let $K$ be obtained by multiplying that block by $\zeta_k$. Then $G_K$ is obtained from $G_H$ by the cyclic switching of Theorem~\ref{thm:cyclic_switching}. In particular, $G_H$ and $G_K$ are cospectral. The balance hypothesis is automatic when $k$ is prime.
\end{corollary}

\begin{proof}
Assume the selected columns are $0,\ldots,m-1$ and put
\[
 X_a=\{C_{j,a}:0\leq j<m\}\qquad(a\in\mathbb Z_k).
\]
Their union is independent. Every selected row vertex $R_{i,a}$, with $0\leq i<q$, is adjacent to all of $X_a$ and none of the other classes. By Lemma \ref{lem:prime_balance}, every other row vertex has $m/k$ neighbours in each $X_a$; other column vertices have none. The graph theorem therefore applies. Its cyclic change replaces the matching for an entry $1$ in the selected block by the matching for $\zeta_k$, giving $G_K$.
\end{proof}

We can also express this correspondence directly in the original vertex order. Put $E=J_m/m$ and define the complex unitary matrix
\[
 L=I_m+(\zeta_k-1)E
\]
and the real orthogonal matrix
\begin{equation}\label{eq:cyclic_lift}
 \widehat L=(I_m-E)\otimes I_k+E\otimes P_k
 =I_{mk}+E\otimes(P_k-I_k).
\end{equation}
The first matrix multiplies a constant row by $\zeta_k$ and fixes a row with sum zero. The second rotates the constant phase-class components and fixes their orthogonal complement. We define $\widehat L$ explicitly, since $\phi_k$ has no additive extension to complex numbers that would justify writing it as $\phi_k(L)$.

Under the balance condition,
\begin{align}
 \phi_k(J_{q,m})\widehat L&=\phi_k(\zeta_kJ_{q,m}),\label{prop:phiJL}\\
 \phi_k(A)\widehat L&=\phi_k(A).\label{prop:phiAL}
\end{align}
For the second identity, each row-block sum of $\phi_k(A)$ is $(m/k)J_k$, which annihilates $P_k-I_k$ on the right. Hence, with $n=qm$,
\[
 \phi_k(K)=\phi_k(H)\bigl(\widehat L\oplus I_{(n-m)k}\bigr).
\]
Writing $\mathcal L=I_{nk}\oplus\widehat L\oplus I_{(n-m)k}$ gives
\[
 \mathcal M_K=\mathcal L^{\top}\mathcal M_H\mathcal L.
\]
This proves the graph correspondence by an explicit real orthogonal similarity. For composite phase, zero row sums alone do not imply the condition of Lemma \ref{lem:prime_balance}; for example, $(1,1,-1,-1)$ has sum zero in $\mu_4$ but contains neither of the phases $\zeta_{4}$, and $-\zeta_{4}$. The balance hypothesis must therefore be retained in the graph lifting statement.

\section{Concluding remarks}\label{section:concluding_remarks}

Our results provide a basis for transferring tools between matrix switching and graph switching. The correspondences in Sections 5 and 6 show how suitable matrix operations lead to graph switchings and how graph-switching methods yield operations on Butson matrices. The spectral analysis in Section 4 also explains a limitation of this transfer: for fixed order and phase, preserving the Butson property guarantees cospectral associated graphs at prime phase, but need not do so at composite phase.

A natural next step is to characterise the graph-switching partitions that preserve the class of Butson Hadamard graphs. Together with structural criteria for spectral preservation at composite phase, such a characterisation could guide the search for inequivalent matrices and nonisomorphic cospectral graphs.

The computational work in this paper, such as graph isomorphism and Hadamard equivalence checking, was carried out with {\sc Magma} \cite{Magma}.

\section*{Acknowledgements}
Aida Abiad acknowledges support from the Dutch Research Council (NWO),
Vidi grant VI.Vidi.213.085 (grant ID \url{https://doi.org/10.61686/CSWNV56135}). This work was done in part while the first author was visiting the Simons Institute for the Theory of Computing at UC Berkeley.

\section*{Declaration of AI use}
ChatGPT (OpenAI, accessed October 2026) assisted with literature checking, proofreading and revision
of the exposition. The authors are responsible for the mathematical content and
the final manuscript.

\bibliographystyle{plain}
\bibliography{references}
\end{document}